\documentclass[11pt]{article}
\usepackage{mathrsfs}
\usepackage{amsthm}
\usepackage{amssymb}
\usepackage{amsmath}
\usepackage{graphicx}
\usepackage{color}
\usepackage{amsfonts}
\usepackage{float}
\usepackage{cite}
\usepackage[text={140mm,210mm},left=35mm,vmarginratio=1:1]{geometry}
\newtheorem{theorem}{Theorem}[section]

\newtheorem{lemma}[theorem]{Lemma}
\newtheorem{proposition}[theorem]{Proposition}

\numberwithin{equation}{section}
\normalsize

\begin{document}
\title{\textbf{Moderate deviations of hitting times of a family of density-dependent Markov chains}}

\author{Yuheng He \thanks{\textbf{E-mail}: heyuheng@bjtu.edu.cn \textbf{Address}: School of Mathematics and Statistics, Beijing Jiaotong University, Beijing 100044, China.}\\ Beijing Jiaotong University\\ Xiaofeng Xue \thanks{\textbf{E-mail}: xfxue@bjtu.edu.cn \textbf{Address}: School of Mathematics and Statistics, Beijing Jiaotong University, Beijing 100044, China.}\\ Beijing Jiaotong University}

\date{}
\maketitle

\noindent {\bf Abstract:} In this paper we are concerned with hitting times of a family of density-dependent Markov chains. A moderate deviation principle of the hitting time is given. The proof of the main theorem relies heavily on moderate deviations of density-dependent Markov chains given in \cite{Xue2021} and upper bounds of large deviations of Markov processes given in \cite{Dupuis1991}. An analogue moderate deviation of the hitting time of the diffusion approximation of the density-dependent Markov chain introduced in \cite{Ethier1986} is also given.

\quad


\section{Introduction and main results}\label{section one}

In this paper we are concerned with hitting times of a family of density-dependent Markov chains. For each integer $n\geq 1$, the density-dependent Markov chain  $\{X_t^n\}_{t\geq 0}$ studied in this paper is a continuous-time Markov process with state space which is a subset of $\mathbb{R}$ and generator $\mathcal{L}_n$ given by
\[
\mathcal{L}_nf(u)=n\sum_{i=1}^Ml_iF_i\left(\frac{u}{n}\right)\left[f(u+l_i)-f(u)\right]
\]
for any $u\in \mathbb{R}$ and bounded $f$ from $\mathbb{R}$ to $\mathbb{R}$, where $M\geq 1$ is a given integer, $l_1, l_2,\ldots, l_M\in \mathbb{R}$
and $F_i$ from $\mathbb{R}$ to $[0, +\infty)$ is a $C^1$ function for $i=1,2,\ldots,M$. That is to say, at any moment $t$, $X_t^n$ jumps to $X_t^n+l_i$ at rate $nF_i\left(\frac{X_t^n}{n}\right)$. Throughout this paper, we assume that
\[
F_i(0)=0 \text{~and~} \sup_{u\in \mathbb{R}}\left|\frac{d}{du}F_i(u)\right|<+\infty
\]
for each $i=1,2,\ldots, M$. We further assume that $X_0^n=nx$ for each $n\geq 1$, where $x\in \mathbb{R}$ is a fixed point.

We recall two important examples of above density-dependent Markov Chains.

\textbf{Example 1} \emph{Birth-and-death processes}. If $X_0^n\in\{0,1,2,\ldots\}$, $M=2$, $l_1=1, l_2=-1$ and $F_1(u)=\lambda u, F_2(u)=\theta u$ for some $\lambda, \theta\geq 0$ and any $u\geq 0$, then our model reduces to the birth-and-death process. This process describes the fluctuation of a population, where each individual independently gives birth to a new one at rate $\lambda$ and dies at rate $\theta$.

\textbf{Example 2} \emph{The susceptible-infected-susceptible epidemic on the complete graph}. If $X_0^n\in \{0,1,\ldots,n\}$, $M=2$, $l_1=1, l_2=-1$ and $F_1(u)=\lambda u(1-u), F_2(u)=\theta u$ for $u\in [0, 1]$, then our model reduces to the susceptible-infected-susceptible (SIS) epidemic model on the complete graph with $n$ vertices. Assuming that a susceptible vertex is infected by an infected one at rate $\frac{\lambda}{n}$ and an infected vertex becomes susceptible at rate $\theta$, then $X_t^n$ describes the number of infected vertices at moment $t$.

Under our assumption, the following proposition is a direct application of the main theorem given in \cite{Kurtz1978}, where the density-dependent Markov chain is first introduced.

\begin{proposition}\label{proposition 1.1}
(Kurtz, 1978, \cite{Kurtz1978}) Let $\{x_t\}_{t\geq 0}$ be the unique solution to the ordinary differential equation
\[
\begin{cases}
&\frac{d}{dt}x_t=\sum_{i=1}^Ml_iF_i(x_t),\\
&x_0=x,
\end{cases}
\]
then
\[
\lim_{n\rightarrow+\infty}\frac{X_t^n}{n}=x_t
\]
in probability for any $t\geq 0$.
\end{proposition}
Without loss of generality, in this paper we assume that $\sum_{i=1}^Ml_iF_i(x)>0$. Then it is reasonable to define
\[
x_\infty=\inf\left\{u>x:~\sum_{i=1}^Ml_iF_i(u)=0\right\}.
\]
Note that we let $x_\infty=+\infty$ if $\left\{u>x:~\sum_{i=1}^Ml_iF_i(u)=0\right\}=\emptyset$. Consequently, $\{x_t\}_{t\geq 0}$ is increasing with $t$ such that $\lim_{t\rightarrow+\infty}x_t=x_\infty$ and there is a unique $\tau_r>0$ such that
\[
x_{\tau_r}=r
\]
for any $r\in (x, x_\infty)$. According to the ordinary differential equation $\{x_t\}_{t\geq 0}$ follows, it is easy to check that
\begin{equation}\label{equ tau}
\tau_r=\displaystyle{\int_x^r\frac{1}{\sum_{i=1}^Ml_iF_i(u)}du}.
\end{equation}

For each $n\geq 1$ and any $r\in (x, x_\infty)$, we define
\[
\tau_r^n=\inf\left\{t:~X_t^n\geq nr\right\}.
\]
Note that $\{X_t^n\}_{t\geq 0}$ is not continuous with $t$, hence the moment when $\frac{X^n}{n}$ `hits' $r$ for the first time should be rigorously defined as the moment when $\frac{X^n}{n}$ exceeds $r$ for the first time.

LLN and CLT of hitting times of density-dependent Markov chains are investigated in Chapter 11 of \cite{Ethier1986}. As a direct application of the main result given in Section 11.4 of \cite{Ethier1986}, we have the following proposition.

\begin{proposition}\label{proposition 1.2 LLN and CLT of hitting times}
(Ethier and Kurtz, 1986, \cite{Ethier1986})

1) For any $r\in (x, x_\infty)$, $\lim_{n\rightarrow+\infty}\tau_r^n=\tau_r$ in probability.

2) For any $r\in (x, x_\infty)$, $\sqrt{n}\left(\tau_r^n-\tau_r\right)$ converges weakly to $-\frac{V_{\tau_r}}{\sum_{i=1}^Ml_iF_i(r)}$ as $n\rightarrow+\infty$, where $\{V_t\}_{t\geq 0}$ satisfies
\[
\begin{cases}
&dV_t=\left(\sum_{i=1}^Ml_iF_i^\prime(x_t)\right)V_tdt+\sqrt{\sum_{i=1}^Ml_i^2F_i(x_t)}d\mathcal{B}_t,\\
&V_0=0,
\end{cases}
\]
where $\{\mathcal{B}_t\}_{t\geq 0}$ is a standard Brownian motion.
\end{proposition}

Large deviations of hitting times of density-dependent Markov chains are also discussed in previously literatures. A complete large deviation principle of hitting times with upper and lower bounds under the assumption that $\log F_i$ are bounded for all $i$ is given in Chapter 5 of \cite{Shwartz1995}. However, in this paper we assume that $F_i(0)=0$, hence we can only give a upper bound of the large deviation of $\tau_r^n$ according to the main theorem given in \cite{Dupuis1991}. For mathematical details, see Section \ref{section Preliminaries}.

According to the above limit theorems of $\tau_r^n$, it is natural to further ask whether $\tau_r^n$ follows some moderate deviation principle. The following main result of this paper gives a positive answer of this question.

\begin{theorem}\label{theorem main moderate deviation}
Let $\{a_n\}_{n\geq 1}$ be a positive sequence that
\[
\lim_{n\rightarrow\infty}\frac{a_n}{n}=\lim_{n\rightarrow+\infty}\frac{\sqrt{n}}{a_n}=0,
\]
then
\begin{equation}\label{equ 1.1}
\lim_{n\rightarrow+\infty}\frac{n}{a_n^2}\log P\left(\frac{n}{a_n}(\tau_r^n-\tau_r)>t\right)=-\frac{t^2}{2\displaystyle{\int_x^r \frac{\sum_{i=1}^Ml_i^2F_i(u)}{\left(\sum_{i=1}^Ml_iF_i(u)\right)^3}du}}
\end{equation}
and
\begin{equation}\label{equ 1.2}
\lim_{n\rightarrow+\infty}\frac{n}{a_n^2}\log P\left(\frac{n}{a_n}(\tau_r^n-\tau_r)<-t\right)=-\frac{t^2}{2\displaystyle{\int_x^r \frac{\sum_{i=1}^Ml_i^2F_i(u)}{\left(\sum_{i=1}^Ml_iF_i(u)\right)^3}du}}
\end{equation}
for any $r\in (x, x_\infty)$ and $t>0$.
\end{theorem}

Section 3 of Chapter 11 of \cite{Ethier1986} gives diffusion approximations $\{Z_t^n\}_{t\geq 0}$ of $\left\{\frac{X_t^n}{n}\right\}_{t\geq 0}$, where
\[
\begin{cases}
&dZ_t^n=\sum_{i=1}^Ml_iF_i(Z_t^n)+\frac{1}{\sqrt{n}}\sqrt{\sum_{i=1}l_i^2F_i(Z_t^n)}d\mathcal{B}_t,\\
& Z_0^n=x.
\end{cases}
\]
It is shown in Section 11.3 of \cite{Ethier1986} that $\{Z_t^n\}_{t\geq 0}$ and $\left\{\frac{X_t^n}{n}\right\}_{t\geq 0}$ can be coupled in a same probability space such that
\[
\lim_{n\rightarrow+\infty}P\left(\sup_{0\leq t\leq T}\left|\frac{X_t^n}{n}-Z_t^n\right|>\epsilon\right)=0
\]
for any $\epsilon>0$ and $T>0$.

Our next main result shows that the hitting time of $\{Z_t^n\}_{t\geq 0}$ follows the same moderate deviation principle as that of $\left\{\frac{X_t^n}{n}\right\}_{t\geq 0}$. In detail, for $r\in (x, x_\infty)$, let
\[
\pi_r^n=\inf\{t>0:~Z_t^n=r\},
\]
then we have the following theorem.

\begin{theorem}\label{theorem main moderate deviation of diffusion approxiamtion}
Let $\{a_n\}_{n\geq 1}$ be a positive sequence that
\[
\lim_{n\rightarrow\infty}\frac{a_n}{n}=\lim_{n\rightarrow+\infty}\frac{\sqrt{n}}{a_n}=0,
\]
then
\begin{equation}\label{equ 1.3}
\lim_{n\rightarrow+\infty}\frac{n}{a_n^2}\log P\left(\frac{n}{a_n}(\pi_r^n-\tau_r)>t\right)=-\frac{t^2}{2\displaystyle{\int_x^r \frac{\sum_{i=1}^Ml_i^2F_i(u)}{\left(\sum_{i=1}^Ml_iF_i(u)\right)^3}du}}
\end{equation}
and
\begin{equation}\label{equ 1.4}
\lim_{n\rightarrow+\infty}\frac{n}{a_n^2}\log P\left(\frac{n}{a_n}(\pi_r^n-\tau_r)<-t\right)=-\frac{t^2}{2\displaystyle{\int_x^r \frac{\sum_{i=1}^Ml_i^2F_i(u)}{\left(\sum_{i=1}^Ml_iF_i(u)\right)^3}du}}
\end{equation}
for any $r\in (x, x_\infty)$ and $t>0$.
\end{theorem}

The proof of Theorem \ref{theorem main moderate deviation} is given in Section \ref{section proof}, which relies heavily on moderate deviations of density-dependent Markov chains given in \cite{Xue2021} and an upper bound of the large deviation of $\tau_r^n$ following from upper bounds of large deviations of discontinuous Markov processes given in \cite{Dupuis1991}. The outline of the proof of Theorem \ref{theorem main moderate deviation of diffusion approxiamtion} is given in Section \ref{section proof outline}, which is an analogue of the proof of Theorem \ref{theorem main moderate deviation} according to large and moderate deviations of diffusion processes given in \cite{Varadhan2016} and \cite{Guillin2003}.

\section{Preliminary results}\label{section Preliminaries}
As a preparation for the proof of Theorem \ref{theorem main moderate deviation}, in this section we recall moderate deviations of density-dependent Markov chains and give a upper bound of the large deviation of $\tau_r^n$.

\subsection{A recall of moderate deviations of density-dependent Markov chains}
For any $t\geq 0$, let $\theta_t^n=\frac{X_t^n-nx_t}{a_n}$, then for any given $T>0$, $\theta^{n,T}:=\{\theta_t^n\}_{0\leq t\leq T}$ is a random element in $\mathcal{D}\left([0, T], \mathbb{R}\right)$, where $\mathcal{D}\left([0, T], \mathbb{R}\right)$ is the set of c\`{a}dl\`{a}g functions from $[0, T]$ to $\mathbb{R}$. According to our assumptions, density-dependent Markov chains in this paper belongs to those investigated in \cite{Xue2021}. Therefore, by Theorem 2.2 of \cite{Xue2021}, we have the following proposition.

\begin{proposition}\label{Proposition 2.1}(Xue, 2021, \cite{Xue2021})
For any $f\in \mathcal{D}\left([0, T], \mathbb{R}\right)$, let
\[
I_T(f)=
\begin{cases}
& \frac{1}{2}\int_0^T\frac{\left(f^\prime(u)-C(u)f(u)\right)^2}{\beta(u)}du \text{\quad if~}f\text{~is absolutely continuous and~}f(0)=0,\\
& +\infty \text{\quad otherwise},
\end{cases}
\]
where $C(u)=\sum_{i=1}^Ml_iF^\prime_i(x_u)$ and $\beta(u)=\sum_{i=1}^Ml_i^2F_i(x_u)$, then
\[
\limsup_{n\rightarrow+\infty}\frac{n}{a_n^2}\log P\left(\theta^{n,T}\in C\right)\leq -\inf_{f\in C}I_T(f)
\]
for any closed set $C\subseteq \mathcal{D}\left([0, T], \mathbb{R}\right)$ and
\[
\liminf_{n\rightarrow+\infty}\frac{n}{a_n^2}\log P\left(\theta^{n, T}\in O\right)\geq -\inf_{f\in O}I_T(f)
\]
for any open set $O\subseteq \mathcal{D}\left([0, T], \mathbb{R}\right)$.
\end{proposition}

The following property of $I_T$ plays key role in the proof of Theorem \ref{theorem main moderate deviation}.

\begin{proposition}\label{proposition 2.2 lower bound of IT}
For any $a\in R$,
\[
\inf\left\{I_T(f):~f\in \mathcal{D}\left([0, T], \mathbb{R}\right)\text{~and~}f(T)=a\right\}=\frac{a^2}{2\int_0^T\beta(u)e^{2\int_u^TC(\rho)d\rho}du}.
\]
\end{proposition}

\proof[Proof of Proposition \ref{proposition 2.2 lower bound of IT}]

For any absolutely continuous $f$ with $f(0)=0$ and $f(T)=a$, let $h_f(u)=\frac{f^\prime(u)-C(u)f(u)}{\sqrt{\beta(u)}}$, then $I_T(f)=\frac{1}{2}\int_0^T h^2_f(u)du$ and
\[
f(t)=\int_0^te^{\int_u^t C(\rho)d\rho}\sqrt{\beta(u)}h_f(u)du
\]
for any $t\geq 0$ by solving the equation $f^\prime(u)-C(u)f(u)=\sqrt{\beta(u)}h_f(u)$. As a result, by Cauchy-Schwarz's inequality,
\begin{align*}
a^2=f^2(T)&=\left(\int_0^Te^{\int_u^T C(\rho)d\rho}\sqrt{\beta(u)}h_f(u)du\right)^2 \\
&\leq \int_0^Te^{2\int_u^T C(\rho)d\rho}\beta(u)du\int_0^Th_f^2(u)du=2I_T(f)\int_0^Te^{2\int_u^T C(\rho)d\rho}\beta(u)du
\end{align*}
and hence
\[
I_T(f)\geq \frac{a^2}{2\int_0^T\beta(u)e^{2\int_u^TC(\rho)d\rho}du}.
\]
On the other hand, let
\[
h_{T, a}(u)=\frac{a\sqrt{\beta(u)}e^{\int_u^TC(s)ds}}{\int_0^T\beta(r)e^{2\int_r^TC(\theta)d\theta}dr}
\]
and $f_{T, a}(t)=\int_0^te^{\int_u^t C(\rho)d\rho}\sqrt{\beta(u)}h_{T,a}(u)du$, then $f_{T, a}(0)=0, f_{T, a}(T)=a$ and
\[
I_T(f_{T,a})=\frac{1}{2}\int_0^T h^2_{T, a}(u)du=\frac{a^2}{2\int_0^T\beta(u)e^{2\int_u^TC(\rho)d\rho}du}
\]
and the proof is complete.

\qed

\subsection{An upper bound of the large deviation of the hitting time}

In this subsection we prove the following lemma, which gives an upper bound of the large deviation of our hitting time.

\begin{lemma}\label{lemma 2.3 LDPupperboundHittingTime}
For any $r\in (x, x_\infty)$ and $\epsilon>0$,
\[
\limsup_{n\rightarrow+\infty}\frac{1}{n}\log P\left(|\tau_r^n-\tau_r|>\epsilon\right)<0.
\]
\end{lemma}

The proof of Lemma \ref{lemma 2.3 LDPupperboundHittingTime} relies on upper bounds of large deviations of a family of discontinuous Markov processes given in \cite{Dupuis1991}. According to our assumptions, density-dependent Markov chains in this paper are examples of stochastic processes investigated in \cite{Dupuis1991}, hence by Theorem 1.1 of \cite{Dupuis1991}, we have the following the proposition.

\begin{proposition}\label{Proposition 2.4 UpperBoundofLargeDeviationofMarkovProceeses}(Dupuis, Ellis and Weiss, 1991, \cite{Dupuis1991})
For any $f\in \mathcal{D}\left([0, T], \mathbb{R}\right)$, let
\[
J_T(f)=
\begin{cases}
& \int_0^Tl(f(u), f^\prime(u))du \text{\quad if~}f\text{~is absolutely continuous and~}f(0)=x,\\
& +\infty \text{\quad otherwise},
\end{cases}
\]
where
\[
l(x,y)=\sup_{b\in \mathbb{R}}\left\{by-\sum_{i=1}^MF_i(x)\left(e^{bl_i}-1\right)\right\},
\]
then $J_T$ is a good rate function and
\[
\limsup_{n\rightarrow+\infty}\frac{1}{n}\log P\left(\left\{\frac{X_t^n}{n}\right\}_{0\leq t\leq T}\in C\right)\leq -\inf_{f\in C}J_T(f)
\]
for any closed $C\subseteq \mathcal{D}\left([0, T], \mathbb{R}\right)$.
\end{proposition}

Note that according to a calculus of variation it is easy to check that $J_T(f)\geq 0$ and $J_T(f)=0$ when and only when $f(t)=x_t$ for any $0\leq t\leq T$.

Now we give the proof of Lemma \ref{lemma 2.3 LDPupperboundHittingTime}.

\proof[Proof of Lemma \ref{lemma 2.3 LDPupperboundHittingTime}]

Without loss of generality, we assume that $0<\epsilon<\tau_r$, then
\[
\left\{\tau_r^n-\tau_r<-\epsilon\right\}\subseteq \left\{\left\{\frac{X_t^n}{n}\right\}_{0\leq t\leq \tau_r-\epsilon}\in A_\epsilon^-\right\},
\]
where
\[
A_\epsilon^-=\left\{f\in \mathcal{D}\left([0, \tau_r-\epsilon], \mathbb{R}\right):~\inf_{0\leq t\leq \tau_r-\epsilon}f_t\geq r\right\}.
\]
It is easy to check that $A_\epsilon^-$ is closed in $\mathcal{D}\left([0, \tau_r-\epsilon], \mathbb{R}\right)$. According to facts that $\{x_t\}_{0\leq t\leq \tau_r-\epsilon}\not\in A_\epsilon^-$, $J_{\tau_r-\epsilon}$ is a good rate function and $J_{\tau_r-\epsilon}(f)>0$ for $\{f_t\}_{0\leq t\leq \tau_r-\epsilon}\neq \{x_t\}_{0\leq t\leq \tau_r-\epsilon}$, we have
\[
\inf_{f\in A_\epsilon^-}J_{\tau_r-\epsilon}(f)>0
\]
and hence
\begin{equation}\label{equ 2.1}
\limsup_{n\rightarrow+\infty}\frac{1}{n}\log P\left(\tau_r^n-\tau_r<-\epsilon\right)\leq -\inf_{f\in A_\epsilon^-}J_{\tau_r-\epsilon}(f)<0
\end{equation}
by Proposition \ref{Proposition 2.4 UpperBoundofLargeDeviationofMarkovProceeses}. Similarly,
\begin{equation*}
\limsup_{n\rightarrow+\infty}\frac{1}{n}\log P\left(\tau_r^n-\tau_r>\epsilon\right)\leq -\inf_{f\in A_\epsilon^+}J_{\tau_r+\epsilon}(f),
\end{equation*}
where
\[
A_\epsilon^+=\left\{f\in \mathcal{D}\left([0, \tau_r+\epsilon], \mathbb{R}\right):~\sup_{0\leq t\leq \tau_r+\epsilon}f_t\leq r\right\}.
\]
It is easy to check that $A_\epsilon^+$ is also closed. Since $r\in (x, x_\infty)$, $x^\prime(\tau_r)>0$ and hence $\{x_t\}_{0\leq t\leq \tau_r+\epsilon}\not\in A_\epsilon^+$. Therefore,
\[
\inf_{f\in A_\epsilon^+}J_{\tau_r+\epsilon}(f)>0
\]
and hence
\begin{equation}\label{equ 2.2}
\limsup_{n\rightarrow+\infty}\frac{1}{n}\log P\left(\tau_r^n-\tau_r>\epsilon\right)<0.
\end{equation}
Lemma \ref{lemma 2.3 LDPupperboundHittingTime} follows from Equations \eqref{equ 2.1} and \eqref{equ 2.2}.

\qed

\section{The proof of Theorem \ref{theorem main moderate deviation}}\label{section proof}

In this section we prove Theorem \ref{theorem main moderate deviation}. Here we only give details of the proof of Equation \eqref{equ 1.1} since Equation \eqref{equ 1.2} can be proved in the same way. Equation \eqref{equ 1.1} follows from
\begin{equation}\label{equ 3.1}
\limsup_{n\rightarrow+\infty}\frac{n}{a_n^2}\log P\left(\frac{n}{a_n}(\tau_r^n-\tau_r)>t\right)\leq-\frac{t^2}{2\displaystyle{\int_x^r \frac{\sum_{i=1}^Ml_i^2F_i(u)}{\left(\sum_{i=1}^Ml_iF_i(u)\right)^3}du}}
\end{equation}
and
\begin{equation}\label{equ 3.2}
\liminf_{n\rightarrow+\infty}\frac{n}{a_n^2}\log P\left(\frac{n}{a_n}(\tau_r^n-\tau_r)>t\right)\geq-\frac{t^2}{2\displaystyle{\int_x^r \frac{\sum_{i=1}^Ml_i^2F_i(u)}{\left(\sum_{i=1}^Ml_iF_i(u)\right)^3}du}}.
\end{equation}
We first prove Equation \eqref{equ 3.1}.

\proof[Proof of Equation \eqref{equ 3.1}]

Since $r\in (x, x_\infty)$, we have $x^\prime(\tau_r)>0$. Then for any $0<\epsilon<x^\prime(\tau_r)$, there exists $\delta_1\in (0, \tau_r)$ such that
\[
|x^{\prime}(s)-x^{\prime}(\tau_r)|<\epsilon
\]
when $|s-\tau_r|\leq \delta_1$. Hence, by Lagrange's mean value theorem,
\[
P\left(\frac{n}{a_n}(\tau_r^n-\tau_r)>t, |\tau_r^n-\tau_r|\leq \delta_1\right)
\leq P\left(\frac{n(x_{\tau^n_r}-x_{\tau_r})}{a_n}>(x^\prime(\tau_r)-\epsilon)t\right).
\]
Then, since $\frac{a_n^2}{n^2}\rightarrow 0$,
\begin{align}\label{equ 3.3}
&\limsup_{n\rightarrow+\infty}\frac{n}{a_n^2}\log P\left(\frac{n}{a_n}(\tau_r^n-\tau_r)>t\right) \\
&\leq \limsup_{n\rightarrow+\infty}\frac{n}{a_n^2}\log P\left(\frac{n(x_{\tau^n_r}-x_{\tau_r})}{a_n}>(x^\prime(\tau_r)-\epsilon)t\right) \notag
\end{align}
according to Lemma \ref{lemma 2.3 LDPupperboundHittingTime}. According to definitions of $\tau_r$ and $\tau_r^n$,
\[
x_{\tau_r}=r \text{~and~}\left|\frac{X_{\tau_r^n}^n}{n}-r\right|\leq \frac{\max\{|l_1|,\ldots, |l_M|\}}{n}
\]
and hence
\[
\left|\frac{n}{a_n}\left(x_{\tau^n_r}-x_{\tau_r}\right)+\frac{X_{\tau_r^n}-nx_{\tau_r^n}}{a_n}\right|\leq \frac{\max\{|l_1|, \ldots, |l_M|\}}{a_n}.
\]
As a result, for any $0<\delta_2<(x^\prime(\tau_r)-\epsilon)t$,
\[
P\left(\frac{n(x_{\tau^n_r}-x_{\tau_r})}{a_n}>(x^\prime(\tau_r)-\epsilon)t\right)\leq P\left(\frac{X^n_{\tau_r^n}-nx_{\tau_r^n}}{a_n}< -\left((x^\prime(\tau_r)-\epsilon)t-\delta_2\right)\right)
\]
when $n$ is sufficiently large. For simplicity, we write $(x^\prime(\tau_r)-\epsilon)t-\delta_2$ as $\zeta(\epsilon, \delta_2)$. By Equation \eqref{equ 3.3},
\begin{align}\label{equ 3.4}
&\limsup_{n\rightarrow+\infty}\frac{n}{a_n^2}\log P\left(\frac{n}{a_n}(\tau_r^n-\tau_r)>t\right) \\
&\leq \limsup_{n\rightarrow+\infty}\frac{n}{a_n^2}\log P\left(\frac{X^n_{\tau_r^n}-nx_{\tau_r^n}}{a_n}<-\zeta(\epsilon, \delta_2)\right). \notag
\end{align}
For any $0<\delta_3<\frac{1}{3}\tau_r$,
\[
\left\{\frac{X^n_{\tau_r^n}-nx_{\tau_r^n}}{a_n}<-\zeta(\epsilon, \delta_2), |\tau_r^n-\tau_r|\leq \delta_3\right\}
\subseteq \left\{\theta^{n,\tau_r+\delta_3}\in H_{\epsilon, \delta_2, \delta_3}\right\},
\]
where $\theta^{n, \tau_r+\delta_3}=\{\frac{X_t^n-nx_t}{a_n}\}_{0\leq t\leq \tau_r+\delta_3}$ as defined in Section \ref{section Preliminaries} and
\[
H_{\epsilon, \delta_2, \delta_3}=\left\{f\in \mathcal{D}\left([0, \tau_r+\delta_3], \mathbb{R}\right):~\inf_{\tau_r-\delta_3\leq s\leq \tau_r+\delta_3}f(s)\leq -\zeta(\epsilon, \delta_2)\right\}.
\]
Let $\overline{H}_{\epsilon, \delta_2, \delta_3}$ be the closure of $H_{\epsilon, \delta_2, \delta_3}$ in $\mathcal{D}\left([0, \tau_r+\delta_3], \mathbb{R}\right)$, then by Equation \eqref{equ 3.4}, Lemma \ref{lemma 2.3 LDPupperboundHittingTime} and Proposition \ref{Proposition 2.1},
\begin{equation}\label{equ 3.5}
\limsup_{n\rightarrow+\infty}\frac{n}{a_n^2}\log P\left(\frac{n}{a_n}(\tau_r^n-\tau_r)>t\right)\leq -\inf_{f\in \overline{H}_{\epsilon, \delta_2, \delta_3}}I_{\tau_r+\delta_3}(f).
\end{equation}
For any given $f\in \overline{H}_{\epsilon, \delta_2, \delta_3}$, there exists a sequence in $H_{\epsilon, \delta_2, \delta_3}$ which converges to $f$ in $\mathcal{D}\left([0, \tau_r+\delta_3], \mathbb{R}\right)$, hence
\[
\inf_{\tau_r-2\delta_3\leq s\leq \tau_r+\delta}f(s)\leq -\zeta(\epsilon, \delta_2)<0.
\]
Hence there exists $s_0\in (\tau_r-2\delta_3, \tau_r+\delta_3)$ such that $f(s_0)\leq -\zeta(\epsilon, \delta_2)<0$. According to the definition of $I_T$ and Proposition \ref{proposition 2.2 lower bound of IT},
\begin{align*}
I_{\tau_r+\delta_3}(f)&\geq I_{s_0}\left(\{f(u)\}_{0\leq u\leq s_0}\right)\\
&\geq \frac{f^2(s_0)}{2\int_0^{s_0}\beta(u)e^{2\int_u^{s_0}C(\rho)d\rho}du}\geq \frac{\zeta^2(\epsilon, \delta_2)}{2\int_0^{s_0}\beta(u)e^{2\int_u^{s_0}C(\rho)d\rho}du}\\
&\geq \frac{\zeta^2(\epsilon, \delta_2)}{2\sup_{\tau_r-2\delta_3\leq s\leq \tau_r+\delta_3}\int_0^{s}\beta(u)e^{2\int_u^{s}C(\rho)d\rho}du}.
\end{align*}
As a result,
\[
\inf_{f\in \overline{H}_{\epsilon, \delta_2, \delta_3}}I_{\tau_r+\delta_3}(f)\geq \frac{\zeta^2(\epsilon, \delta_2)}{2\sup_{\tau_r-2\delta_3\leq s\leq \tau_r+\delta_3}\int_0^{s}\beta(u)e^{2\int_u^{s}C(\rho)d\rho}du}
\]
and hence
\[
\limsup_{n\rightarrow+\infty}\frac{n}{a_n^2}\log P\left(\frac{n}{a_n}(\tau_r^n-\tau_r)>t\right)\leq
-\frac{\zeta^2(\epsilon, \delta_2)}{2\sup_{\tau_r-2\delta_3\leq s\leq \tau_r+\delta_3}\int_0^{s}\beta(u)e^{2\int_u^{s}C(\rho)d\rho}du}
\]
by Equation \eqref{equ 3.5}. Since $\delta_3, \delta_2, \epsilon$ are arbitrary, let them converge to $0$, then
\begin{equation}\label{equ 3.6}
\limsup_{n\rightarrow+\infty}\frac{n}{a_n^2}\log P\left(\frac{n}{a_n}(\tau_r^n-\tau_r)>t\right)\leq
\frac{(x^\prime(\tau_r))^2t^2}{2\int_0^{\tau_r}\beta(u)e^{2\int_u^{\tau_r}C(\rho)d\rho}du}.
\end{equation}
According to the definition of $\beta(u)$ and $C(\rho)$,
\begin{align*}
\int_0^{\tau_r}\beta(u)e^{2\int_u^{\tau_r}C(\rho)d\rho}du=\int_{x}^r\sum_{i=1}^Ml_i^2F_i(s)e^{2\int_{x^{-1}(s)}^{x^{-1}(r)}C(\rho)d\rho}\frac{1}{\sum_{i=1}^Ml_iF_i(s)}ds
\end{align*}
and
\begin{align*}
\int_{x^{-1}(s)}^{x^{-1}(r)}C(\rho)d\rho&=\int_{s}^r\frac{\sum_{i=1}^Ml_iF_i^\prime(v)}{\sum_{i=1}^Ml_iF_i(v)}dv
=\ln\left(\frac{\sum_{i=1}^Ml_iF_i(r)}{\sum_{i=1}^Ml_iF_i(s)}\right)=\ln\left(\frac{x^\prime(\tau_r)}{\sum_{i=1}^Ml_iF_i(s)}\right).
\end{align*}
As a result,
\begin{equation}\label{equ 3.7}
\int_0^{\tau_r}\beta(u)e^{2\int_u^{\tau_r}C(\rho)d\rho}du
=\left(x^\prime(\tau_r)\right)^2\int_x^r\frac{\sum_{i=1}^Ml_i^2F_i(s)}{\left(\sum_{i=1}^Ml_iF_i(s)\right)^3}ds.
\end{equation}
Equation \eqref{equ 3.1} follows from Equations \eqref{equ 3.6} and \eqref{equ 3.7}.

\qed

At last, we prove Equation \eqref{equ 3.2}.

\proof[Proof of Equation \eqref{equ 3.2}]

For any $0<\epsilon<x^\prime(\tau_r)$, let $\delta_1$ be defined as in the proof of Equation \eqref{equ 3.1}, then
\[
\left\{|\tau_r^n-\tau_r|\leq \delta_1\right\}\bigcap \left\{\frac{n}{a_n}\left(x_{\tau_r^n}-x_{\tau_r}\right)> t(x^\prime(\tau_r)+\epsilon)\right\}
\subseteq \left\{\frac{n}{a_n}(\tau_r^n-\tau_r)>t\right\}.
\]
Hence, by Lemma \ref{lemma 2.3 LDPupperboundHittingTime},
\begin{align}\label{equ 3.8}
&\liminf_{n\rightarrow+\infty}\frac{n}{a_n^2}\log P\left(\frac{n}{a_n}(\tau_r^n-\tau_r)>t\right) \\
&\geq \liminf_{n\rightarrow+\infty}\frac{n}{a_n^2}\log P\left(\frac{n(x_{\tau^n_r}-x_{\tau_r})}{a_n}>(x^\prime(\tau_r)+\epsilon)t\right). \notag
\end{align}
Since
\[
\left|\frac{n}{a_n}\left(x_{\tau^n_r}-x_{\tau_r}\right)+\frac{X_{\tau_r^n}-nx_{\tau_r^n}}{a_n}\right|\leq \frac{\max\{|l_1|, \ldots, |l_M|\}}{a_n},
\]
for any $\delta_2>0$, we have
\[
P\left(\frac{n(x_{\tau^n_r}-x_{\tau_r})}{a_n}>(x^\prime(\tau_r)+\epsilon)t\right)\geq P\left(\frac{X^n_{\tau_r^n}-nx_{\tau_r^n}}{a_n}<-\left(x^\prime(\tau_r)+\epsilon\right)t-\delta_2\right)
\]
when $n$ is sufficiently large. Then, by Equation \eqref{equ 3.8},
\begin{align}\label{equ 3.9}
&\liminf_{n\rightarrow+\infty}\frac{n}{a_n^2}\log P\left(\frac{n}{a_n}(\tau_r^n-\tau_r)>t\right) \\
&\geq \liminf_{n\rightarrow+\infty}\frac{n}{a_n^2}\log P\left(\frac{X^n_{\tau_r^n}-nx_{\tau_r^n}}{a_n}<-\left(x^\prime(\tau_r)+\epsilon\right)t-\delta_2\right). \notag
\end{align}
For simplicity, we write $\left(x^\prime(\tau_r)+\epsilon\right)t+\delta_2$ as $\varpi(\epsilon, \delta_2)$. For any $0<\delta_3<\frac{1}{3}\tau_r$,
\[
\left\{\frac{X^n_{\tau_r^n}-nx_{\tau_r^n}}{a_n}<-\varpi(\epsilon, \delta_2), |\tau_r^n-\tau_r|\leq \delta_3\right\}
\supseteq \left\{\theta^{n,\tau_r+\delta_3}\in G_{\epsilon, \delta_2, \delta_3}, |\tau_r^n-\tau_r|\leq \delta_3\right\},
\]
where
\[
G_{\epsilon, \delta_2, \delta_3}=\left\{f\in \mathcal{D}\left([0, \tau_r+\delta_3], \mathbb{R}\right):~\sup_{\tau_r-\delta_3\leq s\leq \tau_r+\delta_3}f(s)<-\varpi(\epsilon, \delta_2)\right\}.
\]
Then, by Equation \eqref{equ 3.9}, Lemma \ref{lemma 2.3 LDPupperboundHittingTime} and Proposition \ref{Proposition 2.1},
\begin{equation}\label{equ 3.10}
\liminf_{n\rightarrow+\infty}\frac{n}{a_n^2}\log P\left(\frac{n}{a_n}(\tau_r^n-\tau_r)>t\right)\geq -\inf_{f\in G^o_{\epsilon, \delta_2, \delta_3}}I_{\tau_r+\delta_3}(f),
\end{equation}
where $G^o_{\epsilon, \delta_2, \delta_3}$ is the interior of $G_{\epsilon, \delta_2, \delta_3}$. For any $T>0, a\in \mathbb{R}$, let $h_{T, a}, f_{T, a}$ be defined as in the proof of Proposition \ref{proposition 2.2 lower bound of IT}, then $h_{T, a}=ah_{T, 1}$. Here we choose $T=\tau_r+\delta_3$ and
\[
a=-\frac{(1+\epsilon)\varpi(\epsilon, \delta_2)}{\inf_{\tau_r-2\delta_3\leq s\leq \tau_r+\delta_3}\int_0^se^{\int_u^sC(\rho)d\rho}\sqrt{\beta(u)}h_{\tau+\delta_3, 1}(u)du},
\]
then
\[
\sup_{\tau_r-2\delta_3\leq s\leq \tau_r+\delta_3}f_{\tau_r+\delta_3, a}(u)\leq -(1+\epsilon)\varpi(\epsilon, \delta_2)<-\varpi(\epsilon, \delta_2)
\]
and hence $f_{\tau_r+\delta_3, a}\in G^o_{\epsilon, \delta_2, \delta_3}$. According to the proof of Proposition \ref{proposition 2.2 lower bound of IT},
\[
I_{\tau_r+\delta_3}(f_{\tau_r+\delta_3})=\frac{a^2}{2\int_0^{\tau_r+\delta_3}\beta(u)e^{2\int_u^{\tau_r+\delta_3}C(\rho)d\rho}du}.
\]
Then, by Equation \eqref{equ 3.10},
\[
\liminf_{n\rightarrow+\infty}\frac{n}{a_n^2}\log P\left(\frac{n}{a_n}(\tau_r^n-\tau_r)>t\right)\geq
-\frac{a^2}{2\int_0^{\tau_r+\delta_3}\beta(u)e^{2\int_u^{\tau_r+\delta_3}C(\rho)d\rho}du}.
\]
Since $\epsilon, \delta_2, \delta_3$ are arbitrary, let them converge to $0$, then
\[
\liminf_{n\rightarrow+\infty}\frac{n}{a_n^2}\log P\left(\frac{n}{a_n}(\tau_r^n-\tau_r)>t\right)\geq -\frac{(x^\prime(\tau_r))^2t^2}{2\int_0^{\tau_r}\beta(u)e^{2\int_u^{\tau_r}C(\rho)d\rho}du}.
\]
As we have shown in the proof of Equation \eqref{equ 3.1},
\[
\int_0^{\tau_r}\beta(u)e^{2\int_u^{\tau_r}C(\rho)d\rho}du=(x^\prime(\tau_r))^2\displaystyle{\int_x^r \frac{\sum_{i=1}^Ml_i^2F_i(u)}{\left(\sum_{i=1}^Ml_iF_i(u)\right)^3}du}
\]
and hence Equation \eqref{equ 3.2} holds.

\qed

\section{Outline of the proof of Theorem \ref{theorem main moderate deviation of diffusion approxiamtion}}\label{section proof outline}

In this section we give the outline of the proof of Theorem \ref{theorem main moderate deviation of diffusion approxiamtion}, which is an analogue of the proof of Theorem \ref{theorem main moderate deviation}.

We first recall large and moderate deviations of diffusion processes introduced in \cite{Varadhan2016} and \cite{Guillin2003}. According to our assumptions, $\{Z_t^n\}_{t\geq 0}$ belongs to diffusion processes investigated in \cite{Varadhan2016} and \cite{Guillin2003}. Hence, by Theorem 3.12 of \cite{Varadhan2016} and Theorem 1 of \cite{Guillin2003}, we have the following two propositions.

\begin{proposition}\label{proposition LDPdiffusion}(Theorem 3.12 of \cite{Varadhan2016})
For any $f\in \mathcal{D}\left([0, T], \mathbb{R}\right)$, let
\[
K_T(f)=
\begin{cases}
& \displaystyle{\int_0^T\frac{\left(f^\prime(u)-\sum_{i=1}^Ml_iF_i\left(f(u)\right)\right)^2}{\sum_{i=1}^Ml_i^2F_i\left(f(u)\right)}du} \\
&\text{\quad\quad if~}f\text{~is absolutely continuous and~}f(0)=x,\\
& +\infty \text{\quad otherwise},
\end{cases}
\]
then $K_T$ is a good rate function and
\[
\limsup_{n\rightarrow+\infty}\frac{1}{n}\log P\left(\left\{Z_t^n\right\}_{0\leq t\leq T}\in C\right)\leq -\inf_{f\in C}K_T(f)
\]
for any closed $C\subseteq \mathcal{D}\left([0, T], \mathbb{R}\right)$ and
\[
\liminf_{n\rightarrow+\infty}\frac{1}{n}\log P\left(\left\{Z_t^n\right\}_{0\leq t\leq T}\in O\right)\geq -\inf_{f\in O}K_T(f)
\]
for any open $O\subseteq \mathcal{D}\left([0, T], \mathbb{R}\right)$.
\end{proposition}

\begin{proposition}\label{propositon MDPdiffusion} (Guillin, 2003, \cite{Guillin2003})
Let $I_T$ be defined as in Section \ref{section Preliminaries}, then
\[
\limsup_{n\rightarrow+\infty}\frac{n}{a_n^2}\log P\left(\left\{\frac{n(Z_t^n-x_t)}{a_n}\right\}_{0\leq t\leq T}\in C\right)\leq -\inf_{f\in C}I_T(f)
\]
for any closed set $C\subseteq \mathcal{D}\left([0, T], \mathbb{R}\right)$ and
\[
\liminf_{n\rightarrow+\infty}\frac{n}{a_n^2}\log P\left(\left\{\frac{n(Z_t^n-x_t)}{a_n}\right\}_{0\leq t\leq T}\in O\right)\geq -\inf_{f\in O}I_T(f)
\]
for any open set $O\subseteq \mathcal{D}\left([0, T], \mathbb{R}\right)$.
\end{proposition}

Now we give the proof of Theorem \ref{theorem main moderate deviation of diffusion approxiamtion}.

\proof[Outline of the proof of Theorem \ref{theorem main moderate deviation of diffusion approxiamtion}]

It is easy to check that $K_T(f)\geq 0$ and
\[
K_T(f)=0 \text{~when and only when~} \{f(t)\}_{0\leq t\leq T}=\{x_t\}_{0\leq t\leq T}.
\]
Hence, according to Proposition \ref{proposition LDPdiffusion} and an analysis similar with that in the proof of Lemma \ref{lemma 2.3 LDPupperboundHittingTime}, we have
\begin{equation}\label{equ analogue of Lemma 2.3}
\limsup_{n\rightarrow+\infty}\frac{1}{n}\log P\left(|\pi_r^n-\tau_r|>\epsilon\right)<0
\end{equation}
for any $\epsilon>0$. According to the definition of $\pi_r^n$ and $\tau_r$,
\begin{equation}\label{equ 4.2}
\frac{n(Z_{\pi_r^n}-x_{\pi_r^n})}{a_n}=-\frac{n(x_{\pi_r^n}-x_{\tau_r})}{a_n},
\end{equation}
which is an analogue of
\[
\left|\frac{n}{a_n}\left(x_{\tau^n_r}-x_{\tau_r}\right)+\frac{X_{\tau_r^n}-nx_{\tau_r^n}}{a_n}\right|\leq \frac{\max\{|l_1|, \ldots, |l_M|\}}{a_n}.
\]
Furthermore, Equation \eqref{equ analogue of Lemma 2.3} is an analogue of Lemma \ref{lemma 2.3 LDPupperboundHittingTime} and Proposition \ref{propositon MDPdiffusion} is an analogue of Proposition \ref{Proposition 2.1}. Hence Theorem \ref{theorem main moderate deviation of diffusion approxiamtion} holds according to an analogue of the analysis given in the proof of Theorem \ref{theorem main moderate deviation}.

\qed

\section{Applications}\label{section applications}

In this section we apply Theorem \ref{theorem main moderate deviation} in Examples 1 and 2 given Section \ref{section one}. Throughout this section we assume that $\{a_n\}_{n\geq 1}$ is a positive sequence that $\lim_{n\rightarrow+\infty}\frac{a_n}{n}=\lim_{n\rightarrow+\infty}\frac{\sqrt{n}}{a_n}=0$.

\textbf{Example 1} \emph{Birth-and-death process}. Let $\lambda>\theta$ and $x=1$, i.e., $X_0^n=n$ for each $n\geq 1$, then by Theorem \ref{theorem main moderate deviation} and Equation \eqref{equ tau},
\[
\lim_{n\rightarrow+\infty}\frac{n}{a_n^2}\log P\left(\left|\tau_r^n-\tau_r\right|>\frac{a_n t}{n}\right)=-I(r, t)
\]
for any $r>1$ and $t>0$, where
\[
\tau_r=\int_1^r\frac{1}{\lambda u} du=\frac{\log r}{\lambda-\theta}
\]
and
\[
I(r, t)=\frac{t^2}{2\int_1^r \frac{(\lambda+\theta)u}{(\lambda-\theta)^3u^3}du}=\frac{(\lambda-\theta)^3t^2}{2(\lambda+\theta)(1-\frac{1}{r})}.
\]

\begin{figure}[H]
\centering
\includegraphics[scale=0.4]{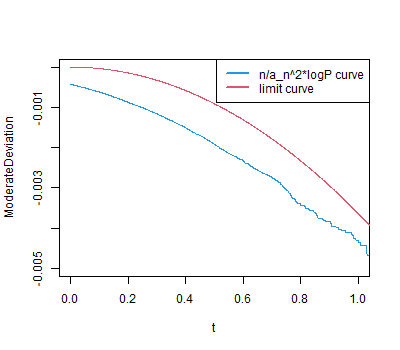}
\caption{$\lambda=1.1$, $\theta=1, x=0.5$ and $r=0.6$}\label{figure 4.1}
\end{figure}

Figure \ref{figure 4.1} gives simulation results under the setting where $\lambda=1.1$, $\theta=1$, $x=1$, $n=10000$, $a_n=n^{0.9}$ and $r=2$. In detail, let $\{X_t^{10000, m}:~t\geq 0\}_{1\leq m\leq 10000}$ be $10000$ independent copies of $\{X_t^{10000}\}_{t\geq 0}$ generated by a computer and $\tau_r^{10000,m}$ be the $\tau_r^{10000}$ of the $m$th copy, then the blue curve in Figure \ref{figure 4.1} gives results of
\[
-\frac{10000}{a_{10000}^2}\log\left(\frac{\sum_{i=1}^{10000}1_{\{\tau_r^{10000,m}-\tau_r>t\}}}{10000}\right)
\]
and the red curve gives the graph of our rate function $-I(2 ,t)$.

\textbf{Example 2} \emph{SIS epidemics on complete graphs}. Let $\theta=1, \lambda>2$, $n$ be even, $x=\frac{1}{2}$ and $\frac{1}{2}<r<\frac{\lambda-1}{\lambda}$, then by Theorem \ref{theorem main moderate deviation} and Equation \eqref{equ tau},
\[
\lim_{n\rightarrow+\infty, \\ \atop n\text{~is even}}\frac{n}{a_n^2}\log P\left(\left|\tau_r^n-\tau_r\right|>\frac{a_n t}{n}\right)=-J(r, t)
\]
for any $t>0$, where
\[
\tau_r=\int_\frac{1}{2}^r\frac{1}{\lambda u(1-u)-u}du=\frac{1}{\lambda-1}\ln(2r)-\frac{1}{\lambda-1}\ln\left(\frac{(\lambda-1)-\lambda r}{\frac{\lambda}{2}-1}\right)
\]
and
\begin{align*}
J(r, t)=\frac{t^2}{2\displaystyle{\int_{\frac{1}{2}}^r \frac{\lambda(1-u)+1}{u^2\left(\lambda(1-u)-1\right)^2}du}}=\frac{\lambda^2t^2}{2\Xi(r)},
\end{align*}
where
\begin{align*}
\Xi(r)=&\frac{\lambda(\lambda+3)}{(\lambda-1)^3}\ln(2r)+\frac{\lambda+1}{(\lambda-1)^2}(2-\frac{1}{r})
+\frac{\lambda(\lambda+3)}{(\lambda-1)^3}\ln\left(\frac{\lambda-1-\frac{1}{2}\lambda}{\lambda-1-\lambda r}\right)\\
&+\frac{2\lambda}{(\lambda-1)^2}\left(\frac{1}{(\lambda-1)-\lambda r}-\frac{1}{(\lambda-1)-\frac{1}{2}r}\right).
\end{align*}

\begin{figure}[H]
\centering
\includegraphics[scale=0.4]{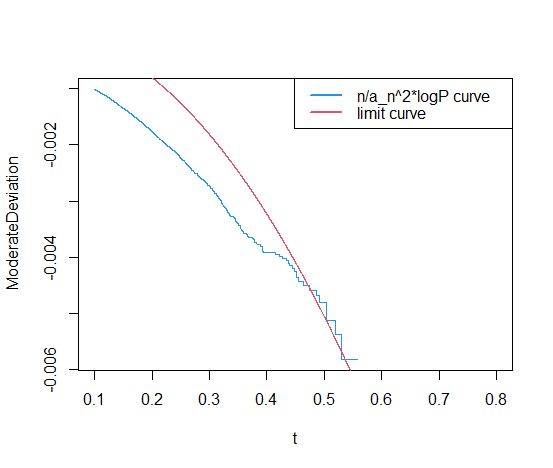}
\caption{$\lambda=3$, $\theta=1, x=0.5$ and $r=0.6$}\label{figure 4.2}
\end{figure}

Figure \ref{figure 4.2} gives simulation results under the setting where $\lambda=3$, $\theta=1$, $x=0.5$, $n=10000$, $a_n=n^{0.9}$ and $r=0.6$. In detail, let $\{X_t^{10000, m}:~t\geq 0\}_{1\leq m\leq 10000}$ be $10000$ independent copies of $\{X_t^{10000}\}_{t\geq 0}$ generated by a computer and $\tau_r^{10000,m}$ be the $\tau_r^{10000}$ of the $m$th copy, then the blue curve in Figure \ref{figure 4.2} gives results of
\[
-\frac{10000}{a_{10000}^2}\log\left(\frac{\sum_{i=1}^{10000}1_{\{\tau_r^{10000,m}-\tau_r>t\}}}{10000}\right)
\]
and the red curve gives the graph of our rate function $-J(0.6 ,t)$.

\quad

\textbf{Acknowledgments.} The authors are grateful to the reviewers. Their comments are great help for the improvement of this paper. The authors are grateful to financial supports from Beijing Jiaotong University with grant number 2022JBMC039 and National Natural Science Foundation of China with
grant number 11501542.

{}
\end{document}